\documentclass[12pt,notitlepage,english]{article}
\pdfoutput=1
\usepackage[utf8]{inputenc}
\usepackage{babel,csquotes}
\usepackage[
backend=biber,
sortcites,style=alphabetic
]{biblatex}
\usepackage{float}
\usepackage[section]{placeins}
\usepackage{graphicx}
\usepackage{mathtools}
\usepackage{amsthm}
\usepackage{amsmath}	
\usepackage{amssymb}
\usepackage{xfrac}
\usepackage{enumitem}
\usepackage{hyperref}
\usepackage{sectsty}
\newcommand{\defeq}{\overset{\textup{{def}}}{=}}

\newcommand{\generate}[1]{\langle #1 \rangle}

\newcommand{\RR}{\mathbb{R}}

\newcommand{\func}[3]{#1: #2 \rightarrow #3}

\newcommand{\eps}{\varepsilon}

\newcommand{\inv}[1]{#1^{-1}}

\newcommand{\st}{:}
\newcommand{\eqdef}{\defeq}
\newcommand{\id}{{1_\Gamma}}
\newcommand{\refl}[1]{\text{Ref}(#1)}
\newcommand{\rot}[1]{\text{Rot}(#1)}
\DeclareMathSymbol{\mhyphen}{\mathord}{AMSa}{"39}
\newtheorem{theorem}{Theorem}[section]
\newtheorem{theorem*}{Theorem}
\newtheorem{proposition}[theorem]{Proposition}
\newtheorem{lemma}[theorem]{Lemma}
\newtheorem{corollary}[theorem]{Corollary}

\theoremstyle{definition}
\newtheorem{definition}[theorem]{Definition}

\theoremstyle{remark}

\newtheorem{remark}[theorem]{Remark}
\title {Illumination in Rational Billiards}
\author{Amit Wolecki}
\sectionfont{\large}
\date{}
\begin{document}
\maketitle
\begin{abstract}
We show that for a rational polygonal billiard, the set of pairs of points that do not illuminate each other (not connected by a billiard trajectory) is finite, and use the same method to extend the results of \cite{EII} and \cite{AW17} about the number of pairs of points that are finitely blocked from each other with a certain blocking cardinality. We rely on previous work about the blocking property in translation surfaces which ultimately stems from results of Eskin, Mirzakhani and Mohammadi on dynamics of moduli spaces of translation surfaces.
\end{abstract}
\section{Introduction}
Consider a source of light as a point in a bounded planar region. The emanating rays hit the boundary and reflect with angle of reflection that equals the angle of incidence, angles taken between the rays and the tangent to the boundary at the point of incidence. 
Does such a light source illuminate the whole region? Is the region illuminable from any point? These two questions of planar geometry are attributed to Ernst Straus in the 1950s, although the requirement that the region is polygonal might have appeared later. An example for a region that is not illuminable by any point was given by Penrose \cite{Pen} in 1958 for a region with curved boundary. 
A set of examples for polygonal regions with two points that do not illuminate each other was given in the 1990s  by Tokarsky \cite{Tok95}. Those polygons have the property of being \emph{rational} in the following sense:
\begin{definition}
A polygon is called rational if all angles between edges are in $\mathbb{Q}\pi$ where the angle between two edges is the angle at the intersection point of their linear continuation.
\end{definition} 
We use the language of billiard dynamical systems to discuss illumination. A pair of points on a polygon $Q$ \emph{do not illuminate each other} if there is no billiard trajectory on $Q$ connecting those points. We consider trajectories between interior points or edge points that are not vertices and assume that billiard trajectories stop upon hitting a vertex. When $Q$ is rational, the process (which will be discussed in detail below) of unfolding $Q$ results in an associated translation surface $M$. Illumination in rational polygons is then generalized by a restrictive condition defined for billiard trajectories as well as translation surface geodesics, called the \emph{blocking property}. A pair of points $A,B \in Q$ is \emph{finitely blocked} if all billiard trajectories between them pass through some finite set, called a \emph{blocking set}. The blocking set is not allowed to contain $A$ and $B$. We denote the minimal cardinality of a blocking set for $A,B$ as $BC(A,B)$. Similarly, blocking and illumination is defined on $M$: a pair $x,y \in M$ of non-singular points is \emph{finitely blocked} if all geodesic trajectories on $M$ connecting $x$ and $y$ pass through some finite set disjoint from $\{x,y\}$. Here we assume that geodesics are not continued through singular points, and discuss continuation through removable singular points later on. The terms blocking set and blocking cardinality are defined similarly, and $BC(x,y)=0$ amounts to non-illumination of $x$ and $y$. The illuminability or blocking of vertices in the billiard or singularities of the surface are not considered in the main results of this paper. However, as the set of pairs of such is finite, including vertices or singularities in the definition of illumination or blocking does not change the main results presented here. 
Also note that illumination and blocking are properties of translation surfaces in general, regardless of whether they are an unfolding of a rational billiard.
 Leli\`evre, Monteil and Weiss \cite{EII} had several results regarding translation surfaces in general: they showed that a translation surface $M$ is a ramified translation cover of the torus if and only if all pairs of non-singular points are finitely blocked from each other and in that case there is a typical $n$ such that every pair of points has $BC\leq n$. They also showed that for a fixed $x \in M$, the set of points that do not illuminate $x$ is finite. For a rational polygon $Q$ that unfolds to a translation surface that is not a torus cover, Apisa and Wright \cite{AW17} showed that if all angles are multiples of $\pi/2$ then each point in $Q$ is finitely blocked from only finitely many other points, and in case there is an angle that is not a multiple of $\pi/2$, then there are at most finitely many pairs of points that are finitely blocked.
\begin{theorem}
\label{thm:main}
Let $Q$ be a rational polygon with connected interior and let $M$ the translation surface obtained by unfolding $Q$, then:
\begin{enumerate}[label=(\roman*)]
	\item \label{thm:main:1} Only finitely many pairs of points in $Q$ do not illuminate each other.
	\item \label{thm:main:2} $\{(A,B) \in Q^2 \st BC(A,B)\leq n\}$ is finite for every integer $n > 0$ if and only if $M$ is not a torus cover.
\end{enumerate}
\end{theorem}
Theorem \ref{thm:main}\ref{thm:main:2} extends the results of \cite{EII} and \cite{AW17} about the possible cardinality of finitely blocked pairs of points in the non-torus cover case, while the result in \cite{EII} grants the "only if" direction. It is concluded in \cite{EII} that for a rational polygon $Q$ and a fixed $A \in Q$, there are at most finitely many points that do not illuminate $A$. Theorem \ref{thm:main}\ref{thm:main:1} extends this result in that the total amount of pairs of points that do not illuminate each other is finite (regardless of whether the unfolded surface is a torus cover or not). 
The results of \cite{EII} and \cite{AW17} rely on work of Eskin, Mirzakhani and Mohammadi \cite{EMM}, applied to the moduli spaces of translation surfaces to describe sets of interest (i.e, subsets of finitely blocked pairs) as submanifolds that locally are complex linear manifolds defined over $\RR$. 
 
Note that the statements in Theorem \ref{thm:main} do not necessarily hold when substituting $Q$ with $M$. Hubert, Schmoll and Troubetzkoy (\cite{HST} \textsection 8) have constructed an example of a non-torus cover translation surface that has infinitely many pairs of points that do not illuminate each other. For other examples see \cite{EII} \textsection 6.3.  

Central to our discussion is the effect of the action of a dihedral group on a translation surface, defined when the surface is an unfolding of a rational polygon. Our results are thus relevant to translation surfaces arising from polygons, but not to arbitrary translation surfaces.
One of the polygonal examples presented in \cite{Tok95} was notable for being a rational polygon with two non-illuminable points, at the time of its discovery it had the fewest number of edges - $26$, compared to other polygonal examples with larger numbers of edges. Castro \cite{Cas97} later presented an example with $24$ edges (published as a Quantum Magazine reader's feedback, in response to an article about Tokarsky's result on the illumination problem), and as a digestif for this paper we exhibit a more efficient example with $22$ edges and justify the following:  
\begin{theorem}
\label{thm:22min}
The minimal number of edges of a polygon with at least one pair of points that do not illuminate each other is at most $22$.
\end{theorem} 
The mentioned examples for non-illumination rely on Tokarsky's main lemma, which we reprove and extend in Lemma \ref{lem:tokrep} to apply to a wider family of examples.
\subsection*{Acknowledgements}
We thank Paul Apisa for an important input regarding the method used in the proof of Lemma \ref{no-inf-intersect-same}. We were recently informed by Thierry Monteil that he also found an example (unpublished) of a polygon with 22 sides that has a pair of points that do not illuminate each other. 
The referee's effort in the review process is acknowledged, and the author is grateful for their insightful comments that have improved the quality and clarity of this paper.
This work is part of the author's MSc thesis conducted at Tel Aviv University under the supervision of Barak Weiss. The author wishes to thank him for the dedicated mentorship, guidance and support throughout this research journey.
The support of BSF grant 2016256 and ISF grant 2095/15 is gratefully 
acknowledged.
\section{Illumination and Blocking in Rational Billiards and Translation Surfaces} 
Let $\Gamma$ be the group generated by the linear parts of reflections along edges of a polygon $Q$. $\Gamma$ is finite if and only if $Q$ is rational and in this case a process of unfolding (also named Katok-Zemliakov construction) yields a translation surface $M$: a compact, orientable surface with an atlas of planar charts (from open subsets of $M \setminus  \Sigma$ to $\mathbb{C}$) such that all transition maps are translations, and equipped with a flat metric defined outside a finite set of singularities $\Sigma$.
For a general description of translation surfaces we refer to \cite{MaTa}, \cite{Vor96}. Also note that $\Gamma$ can be defined for a general polygon, and finite $\Gamma$ is commonly the definition of a rational polygon \cite{MaTa}, \cite{Vor96}.
The unfolding process begins with taking all reflections and rotations of $Q$ by elements of $\Gamma$ with arbitrary translations to avoid overlaps, marking the resulting set of polygons $\{Q_{(\gamma)} : \gamma \in \Gamma\}$ and denoting the affine maps $\phi_\gamma(x) = \gamma x + v_\gamma$ for suitable $v_\gamma \in \RR^2$ such that $\phi_\gamma (Q) = Q_{(\gamma)}$. We then identify edges in the following way: edges $e_1,e_2$ on $Q_{(\gamma_1)}, Q_{(\gamma_2)}$ respectively are glued by translations whenever (1) they correspond to the same edge in $Q$, that is $\inv \phi_{\gamma_1}(e_1) = \inv \phi_{\gamma_2}(e_2)$ and (2) $e_1$ is parallel to $e_2$ and $\gamma_1 = \tau\gamma_2$ where $\tau$ is the derivative of reflection in $e_1$.
For a more detailed description of this unfolding construction we refer to \cite{Vor96} \textsection 2.2.
The following lemma is a standard claim in billiards and we leave the proof to the reader.
\begin{lemma}
\label{lem:exiff}
Let $Q$ be a rational polygon and let $M$ be its unfolding. Let $A, B \in Q$ a pair of distinct points that are not vertices of $Q$, and $E \subset Q \setminus \{A,B\}$. A billiard trajectory on $Q \setminus E$ connecting $A$ and $B$ exists if and only if there exists a geodesic line on $M \setminus (\inv\pi(E) \cup \Sigma)$ connecting a point in $\inv\pi(A)$ to a point in $\inv\pi(B)$ where $\func{\pi}{M}{Q}$ is the natural projection that maps $x \in M$ to a point in $Q$ that unfolds to $x$.
\begin{remark} (Illumination and removable singularities.)
Traditionally, a billiard trajectory is defined such that it stops upon hitting a vertex. When the billiard has a vertex with interior angle $\pi/n$, the billiard flow can be continued to pass through this vertex. The continuation of a billiard trajectory $\alpha(t)$ through a $\pi/n$ vertex at time $t_0$ amounts to the trajectory bouncing back on itself, i.e. $\alpha(t_0+t) = \alpha(t_0-t)$ for all $t$. A vertex with interior angle $\pi/n$ is lifted to a removable singular point in the unfolded surface, on which the flat metric can be defined. A lifted trajectory $\widetilde{\alpha}(t)$ allowed through the removable singularity at $\widetilde{\alpha}(t_0)$ would pass from one copy of the polygon to the copy shifted by a non-trivial element $\gamma \in \Gamma$ such that $\widetilde{\alpha}(t_0+t) = \gamma \widetilde{\alpha}(t_0-t)$ for all $t$. The assertion in Lemma \ref{lem:exiff} holds for billiard and geodesic flows continued through $\pi/n$ vertices and removable singularities respectively. Yet in the general illumination question there is a slight difference between billiards and surfaces. A pair of points $A, B \in Q$ that do not illuminate each other in the traditional billiard flow would still not illuminate each other if trajectories are allowed to continue through $\pi/n$ vertices. However on the unfolded surface we can consider a case where a pair of points in $M$ illuminate each other only when a singularity is removed, these points must descend to the same point in the billiard table $Q$.
\end{remark}
\end{lemma}
\begin{lemma}
\label{lem:nonex} Let $Q$ be a rational polygon, $A,B \in Q$. The following are equivalent:
\begin{enumerate}
	\item $A$ and $B$ are finitely blocked \label{item:Q}
	\item $x$ and $y$ are finitely blocked for all $x \in \inv\pi(A)$ and $y \in \inv\pi(B)$. \label{item:M}
\end{enumerate}
In addition, for all $n \geq 0, x \in \inv\pi(A), y \in \inv\pi(B)$:
\begin{equation} \label{eq:bcn}
	BC(A,B) \leq n \Rightarrow BC(x , y) \leq |\Gamma|n
\end{equation}
\end{lemma}
\begin{proof}
If $E \subset Q \setminus \{A,B\}$ is a blocking set for $A,B$ then by Lemma \ref{lem:exiff} $\inv \pi(E)$ is a blocking set for every pair $x \in \inv\pi(A)$ and $y \in \inv\pi(B)$ and when $|E| \leq n$ for some finite $n$ then $|\inv \pi(E)| \leq |\Gamma||E| \leq |\Gamma|n$.
Conversely, it is enough to prove that for a fixed $x \in \inv\pi(A)$ there is a set $K$ that is a blocking set for any pair of points $(x,y)$ with $y \in \inv{\pi}(B)$, that is, $K$ blocks any geodesic trajectory between $x$ and $y$, does not contain $x$ and is disjoint from $\inv{\pi}(B)$. The claim will then follow by taking $E = \pi(\Gamma K)$ and observing that it does not contain $A$ or $B$, so $\inv\pi(E)$ satisfies a finite blocking set for all pairs in $\inv\pi(A) \times \inv\pi(B)$, and so by Lemma \ref{lem:exiff} $A$ and $B$ are finitely blocked by $E$.
We now fix $x \in \inv\pi(A)$ and construct $K$ iteratively, ensuring that $K$ is disjoint from $\inv{\pi}(B)$. Recall that a blocking set for $(x,y)$ does not contain $x$ or $y$ by definition, and note that by construction $K$ will not contain $x$ as the construction relies on unions of blocking sets. Fix $\{\gamma_1, \dots, \gamma_{|\Gamma|}\}$ an ordering of elements of $\Gamma$ and denoting $y_i$ the point in $\inv\pi(B) \cap Q_{(\gamma_i)}$ we take $K_i$ to be a blocking set for the pair $(x,y_i)$. Now take $K^{(1)} = K_1$, by definition this is a blocking set for $(x,y_1)$. For $1 < i \leq |\Gamma|$, define $K^{(i)} = (K^{(i-1)} \setminus \{y_i\}) \cup (K_i \setminus \{y_1, \dots, y_{i-1}\})$. By induction, let $1 < i \leq |\Gamma|$ and assume that $K^{(i-1)}$ is a blocking set for all pairs of points in $\{(x,y_1), \dots, (x,y_{i-1})\}$. We show that $K^{(i)}$ is a blocking set for all the pairs $\{(x,y_1), \dots, (x,y_{i})\}$. Clearly, $K^{(i)}$ does not contain any element of $\{y_1, \dots, y_{i-1}\}$ as $K^{i-1}$ is by assumption disjoint form $\{y_1, \dots, y_{i-1}\}$, and $y_i$ by definition is not contained in $K_i$. It is required to show that $K^{(i)}$ is blocking all geodesic trajectories between every pair of points $(x,y_1), \dots, (x,y_i)$. First let $\sigma_j$ a geodesic trajectory from $x$ to $y_j$ for some $1 \leq j < i$, and we later deal with the case of a geodesic trajectory $\sigma_i$ from $x$ to $y_i$.
If $\sigma_j$ does not intersect any element of $\{y_1, \dots, y_i\} \setminus \{y_j\}$ along the way to $y_j$, then it will by assumption intersect the contribution of $K^{(i-1)}$ to the union in $K^{(i)}$. Otherwise, take $1 \leq l \leq i$ such that $y_l \neq y_j$ is the first element of $\{y_1, \dots, y_i\} \setminus \{y_j\}$ that $\sigma_j$ intersects. If $l \neq i$ then again $\sigma_j$ will intersect $K^{(i-1)} \setminus \{y_i\}$ before reaching $y_l$, and in case $l=i$ then it will intersect $(K_i \setminus \{y_1, \dots, y_{i-1}\})$ before reaching $y_i$. So in total, $\sigma_j$ intersects $K^{(i)}$ and we are left to deal with $\sigma_i$ a geodesic trajectory between $x$ and $y_i$. Indeed, if $\sigma_i$ passes through a point in $\{y_1, \dots, y_{i-1}\}$ then $\sigma_i$ is intercepted by $K^{(i-1)} \setminus \{y_i\}$ before it reaches any point in $\{y_1, \dots, y_{i-1}\}$. If $\sigma_i$ does not pass through any $\{y_1, \dots, y_{i-1}\}$ then it must pass through a point in $K_i \setminus \{y_1, \dots, y_{i-1}\}$. This amounts to $K^{(i)}$ a blocking set for every pair in $\{(x,y_1), \dots, (x,y_{i})\}$.
Now by induction this holds for $i=|\Gamma|$ and since $\inv{\pi}(B) = \{y_1, \dots, y_{|\Gamma|}\}$ we can take $K = K^{(|\Gamma|)}$ a blocking set for every pair of $(x,y)$ with $y \in \inv{\pi}(B)$.
\end{proof}
\begin{lemma}
  \label{lem:prod-grop-invariance}
  Let $n \geq 0$, the set of lifted pairs of blocked points $\{(x,y) \in (M \setminus \Sigma)^2 \st BC(\pi(x), \pi(y)) \leq n\}$ is invariant under the action of the product group $\Gamma^2$.
  \begin{proof}
    For $x \in M$ nonsingular and $\gamma \in \Gamma$ we have that $\gamma x$ is nonsingular and $\pi(\gamma x) = \pi(x)$ by construction as $M$ is defined by gluing copies of $Q$ that correspond to transformations by elements of $\Gamma$. Observing the action of the product group, we can say that for every $x, y \in M$ nonsingular points $\gamma_1, \gamma_2 \in \Gamma$ we have $BC(\pi(x), \pi(y)) \leq n \Rightarrow BC(\pi(\gamma_1 x), \pi(\gamma_2 y)) \leq n$ which yields the invariance.
  \end{proof}
\end{lemma}
\begin{definition}
  Let $Q$ a rational polygon with connected interior, and let $M$ the translation surface obtained by unfolding $Q$. We denote the set of pairs of points on $Q$ that are blocked with blocking cardinality $n$ as $\mathcal{Q}_n \eqdef \{ (A,B) \in Q^2 \st BC(A,B) \leq n \}$. The set of pairs of points blocked with blocking cardinality $n$ on the translation surface is denoted $\mathcal{M}_n \defeq \{(x,y) \in \widehat{M}^2 \st BC(x,y) \leq n \}$ where $\widehat{M}^2 \eqdef (M \setminus \Sigma)^2$.
\end{definition}
For a general translation surface $M$ we treat the product $\widehat{M}^2$ as a translation structure with charts in $\mathbb{C}^2$, that are induced naturally by the charts of $M \setminus \Sigma$. It is shown in \cite{EII} that $M$ is a branched cover of the torus if and only if there is some $k > 0$ for which $\mathcal{M}_k = \widehat{M}^2$. More details on $\mathcal{M}_k$ in case it is a proper subset of $\widehat{M}^2$ are provided in \cite{EII} and also recalled in Proposition \ref{prop:EII}, describing it as a finite union of linear submanifolds of $\widehat{M}^2$. Note that the real dimension of the complex submanifolds of $\widehat{M}^2$ is considered.
A $2$-dimensional linear submanifold $S \subset \widehat{M}^2$ can be locally defined by affine equations $ax + by = u$ with $x$ and $y$ in the complex coordinates of $M$, $a,b \in \RR$ and $u \in \mathbb{C}$ corresponding to the translation charts of $S$. If $S$ is connected, the coefficients in these affine equations are constant through the charts. Note that $S$ respects the translation structure of $\widehat{M}^2$.
\begin{definition}
 When $S$ is a connected $2$-dimensional linear submanifold of $\widehat{M}^2$ that forms a compact and orientable surface, locally defined by affine equations of the form $ax + by = u$ as described above, such that both real linear coefficients are nonzero, we call $S$ a \emph{translation surface affinely embedded in $\widehat{M}^2$}.
\end{definition}
\begin{proposition} [Leli\`evre, Monteil, Weiss {\cite{EII}}, Theorem 2 and Theorem 11]
\label{prop:EII}
For a translation surface $M$ and an integer $k \geq 0$, the set $\mathcal{M}_k$ either equals $\widehat M^2$ or is a finite union of $0$-manifolds, and $2$-dimensional linear submanifolds of $\widehat M^2$ of the following kinds:
\begin{enumerate}
	\item $F \times M$ or $M \times F$ for a finite $F \subset M$.
  \item A translation surface affinely embedded in $\widehat{M}^2$.
\end{enumerate}
\end{proposition} 
\begin{remark}
Proposition \ref{prop:EII} holds for blocking in general translation surfaces and the proof in \cite{EII} relies on a result by Eskin, Mirzakhani and Mohammadi describing orbit closures of the $SL(2,\RR)$ action on moduli spaces of translation surfaces. An orbit closure appears as a submanifold of the moduli space with an atlas of charts where each chart's image is an open set in an affine space in which the linear part is a complex vector space defined over $\RR$, and so it has an even dimension. The set of all such manifolds that arise from orbit closures is countable, and finite in the case of such submanifolds arising from marking pairs of points that are also blocked with blocking cardinality $n \geq 0$ (as shown in \cite{EII} and cited in Proposition \ref{prop:EII}).
\end{remark}
\begin{lemma}
  \label{no-inf-intersect-same}
  Let $M$ a translation surface obtained by the unfolding of a rational polygon with connected interior $Q$, $k \geq 0$ an integer and let $\mathcal{M}_{Q,k} \eqdef \{(x,y) \in \widehat{M}^2 \st (\pi(x), \pi(y)) \in \mathcal{Q}_k\}$ the set of pairs of points that are lifts of $\mathcal{Q}_k$. If $\mathcal{M}_{|\Gamma|k} \subsetneq \widehat{M}^2$ then for any $S \subset \mathcal{M}_{|\Gamma|k}$ a translation surface affinely embedded in $\widehat{M}^2$, the intersection of $\mathcal{M}_{Q,k}$ with $S$ is finite.
\end{lemma}
\begin{proof}
First note that by Lemma \ref{lem:nonex}, $\mathcal{M}_{Q,k} \subset \mathcal{M}_{|\Gamma|k}$. Assume $\mathcal{M}_{|\Gamma|k} \subsetneq \widehat{M}^2$. By Proposition \ref{prop:EII} it is a finite union of $0$-manifolds and $2$-dimensional linear submanifolds which are either of the form $F \times M$ or $M \times F$ for some finite $F \subset M$, or translation surfaces affinely embedded in $\widehat{M}^2$. Let $S \subset \mathcal{M}_{|\Gamma|k}$ be a translation surface affinely embedded in $\widehat{M}^2$ and assume by contradiction that $\mathcal{M}_{Q,k} \cap S$ is infinite. In the complex coordinates of $\widehat{M}^2$, a suitable neighborhood $U$ of $S$ is a set of solutions for a linear equation of the form $ax + by = u$ for nonzero $a,b \in \mathbb{R}$, with $x$ a point in the first component and $y$ a point in the second component of $\widehat{M}^2$. $U$ can be taken to be an open ball that contains infinitely many elements of $\mathcal{M}_{Q,k}$. Without loss of generality, $Q$ is such that one of its edges is horizontal (the blocking property as well as $Q$ being a rational polygon would persist if we rotate $Q$ and its unfolding). By the unfolding process, the horizontal reflection $R_h$ is an element of $\Gamma$. Consider the action of $\Gamma^2$ on the product space $\widehat{M}^2$, under which $\mathcal{M}_{Q,k}$ is invariant (Lemma \ref{lem:prod-grop-invariance}). By Lemma \ref{lem:nonex} We have that $(R_h, \id)(\mathcal{M}_{Q,k} \cap U)$ is also a subset of $\mathcal{M}_{|\Gamma|k}$. As both components project to infinitely many points it must embed in a translation surface affinely embedded in $\widehat{M}^2$, denoted $S'$ (as described in Proposition \ref{prop:EII}). This yields another real-linear equation that is satisfied in a suitable open ball $N \subset (R_h, \id)U$ of $S'$, so we have that the planar region $\inv{(R_h,\id)}N \cap U$ satisfies the following equations simultaneously:
  \begin{equation*}
    \begin{cases}
    ax + by = u \\
    a'R_hx + b'y = w
    \end{cases}
  \end{equation*}
  For some $u, w \in \mathbb{C}$ with $a, a', b, b'$ nonzero real coefficients. Considered as a set of real-linear equations in $\mathbb{R}^4$ that are satisfied in the $2$-dimensional sub-domain $\inv{(R_h,\id)}N \cap U$, the rank of following matrix equals $2$: 
  \begin{equation*} 
    rk \begin{pmatrix} a & 0 & b & 0 \\
    0 & a & 0 & b \\
    a' & 0 & b' & 0 \\
    0 & -a' & 0 & b'
   \end{pmatrix} = 2
  \end{equation*}
  This matrix can be transformed to the upper triangular form:
  \begin{equation*} 
  \begin{pmatrix} \quad a \quad\;\; & 0\;\; & b & \!\!\!\!\!\!0 \\[3pt]
    \quad0\quad\;\; & a\;\; & 0 & \!\!\!\!\!\!b \\[3pt]
    \quad0\quad\;\;& 0\;\; & b' - \frac{ba'}{a} & \!\!\!\!\!\!0 \\[3pt]
    \quad0\quad\;\;& 0\;\; & 0 & \!\!\!\!\!\!b' + \frac{ba'}{a}
   \end{pmatrix}
  \end{equation*}
  Where each of the last two rows of this matrix cannot be a linear combination of the other rows, hence must be all zeros. Taking $c = ba'/a$ which is nonzero as $a,b,a'$ are nonzero, we get: 
  \begin{align*}
    b' - c = b' + c \\
    - c = c \\
    -1 = 1
  \end{align*}
\end{proof}

\begin{proposition}[Apisa, Wright \cite{AW17}, Corollary 3.8] \label{prop:AW17}
A point $x \in M \setminus \Sigma$ on a translation surface $M$ that is not a branched cover of the torus, is finitely blocked from only finitely many other points.
\end{proposition} 
\begin{corollary}
\label{cor:nonTcover}
If $M$ is not a torus cover, then for every integer $k \geq 0$, $\mathcal{M}_k$ is the union of a finite set and finitely many translation surfaces affinely embedded in $\widehat{M}^2$.
\begin{proof}
By Proposition \ref{prop:AW17}, for every $x \in M \setminus \Sigma$ the set $\{y \in M \setminus \Sigma \st BC(x,y) \leq k\}$
is finite hence $\mathcal{M}_k$ does not coincide with $\widehat{M}^2$ so by Proposition \ref{prop:EII} it is a finite union of $0$ and $2$-dimensional submanifolds of the kinds described in this proposition. Again, by Proposition \ref{prop:AW17}, $\mathcal{M}_k$ cannot contain any set of the form $\{x\} \times M$ or $M \times \{x\}$ for any $x$, so the only $2$-dimensional manifolds that comprise the union are translation surfaces affinely embedded in $\widehat{M}^2$.
\end{proof}
\end{corollary}
\begin{proof}[Proof of Theorem \ref{thm:main}]
  First we note that for $n \geq 0$ we have by Lemma \ref{lem:nonex} that the lift of $\mathcal{Q}_n$, denoted $\mathcal{M}_{Q,n}$ is a subset of $\mathcal{M}_{|\Gamma|n}$, and that it is sufficient to show that $\mathcal{M}_{Q,n}$ is finite as then $\mathcal{Q}_n$ must be finite as a projection of a finite set.
  $(i)$: In the case of illumination, we can always find a pair of close enough points so that they illuminate each other, hence $\mathcal{M}_0$ cannot coincide with $\widehat{M}^2$. We can then consider $\mathcal{M}_0$ as a finite union of $0$ and $2$-dimensional submanifolds, and observe that $\mathcal{M}_0$ cannot contain any submanifolds of the form $F \times M$ or $M \times F$ for some finite $F \subset M$ as these submanifolds contain pairs of close enough points that illuminate each other. By Lemma \ref{no-inf-intersect-same}, $\mathcal{M}_{Q,n}$ has a finite intersection with each of the finitely many possible submanifolds of $\mathcal{M}_0$ which are translation surfaces affinely embedded in $\widehat{M}^2$, so $\mathcal{M}_{Q,n}$ must be finite and hence $\mathcal{Q}_0$ is finite. In other words, there are only finitely many pairs of points on the polygon $Q$ that do not illuminate each other.
  $(ii)$: let $n > 0$ and assume that $M$ is not a torus cover, then by Proposition \ref{prop:AW17}, $\mathcal{M}_{|\Gamma|n}$ does not coincide with $\widehat{M}^2$ and by Corollary \ref{cor:nonTcover} it is the union of a finite set and finitely many translation surfaces affinely embedded in $\widehat{M}^2$. Now similarly, by Lemma \ref{no-inf-intersect-same} $\mathcal{M}_{Q,n}$ has a finite intersection with any translation surface affinely embedded in $\widehat{M}^2$ and hence it is finite and so is $\mathcal{Q}_n$.
  The "only if" direction is evident from the result by Leli\`evre, Monteil and Weiss, that $M$ is a torus cover if and only if for some $k \geq 0$, $\mathcal{M}_k=\widehat{M}^2$ (\cite{EII}, Theorem 1). Note that this statement in \cite{EII} holds for general translation surfaces, that are not necessarily constructed by unfolding a rational polygon. In our case, condition \ref{item:M} in Lemma \ref{lem:nonex} is satisfied for infinitely many distinct $A,B \in Q$, hence $\mathcal{Q}_n$ is infinite. Note that Lemma \ref{lem:nonex} uses a general construction of a blocking set to derive blocking of points $A,B$ in the polygon from blocking in the preimages of $A$ and $B$ on the surface, while in the special case of ramified covers of the torus, \cite{EII} shows explicit constructions that rely on midpoints of trajectories between a pair of points on the torus, which can be lifted to a blocking set for pairs of points in the cover and provide a more explicit build of blocking sets for $\mathcal{Q}_n$.
\end{proof}
\section{A little more about illumination} 
\begin{figure}[ht] 
\center
\includegraphics[width = 0.5\textwidth]{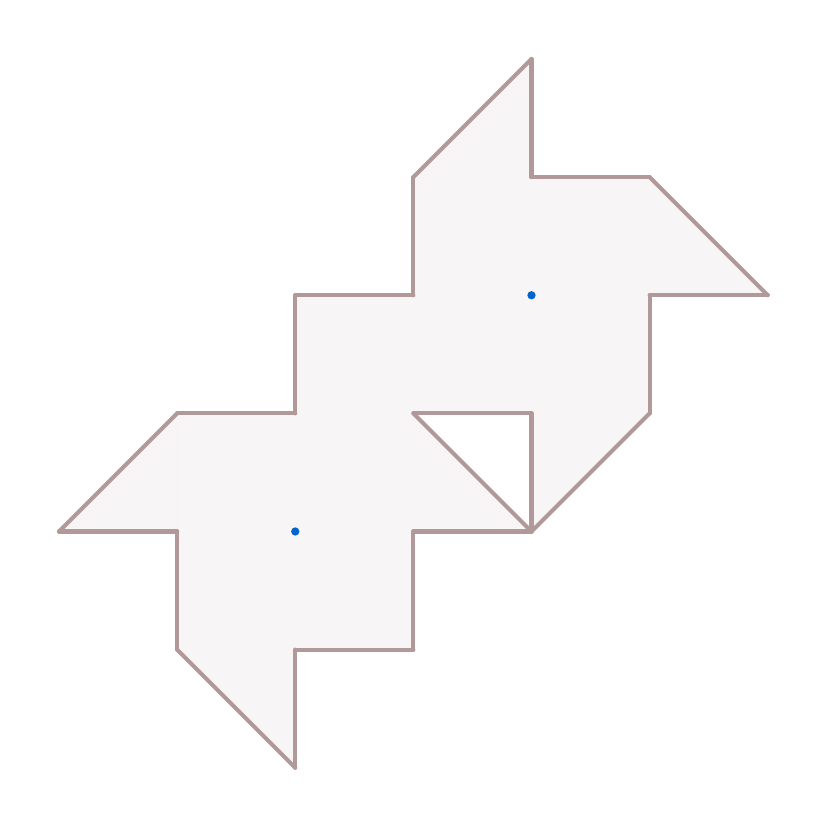}
\caption{A $22$-sided polygon with a pair of points that do not illuminate each other}
\label{fig:22gon}
\end{figure}

The two points on the polygon in Figure \ref{fig:22gon} do not illuminate each other. This polygon's construction is made to work with the same construction theorem that bases the proof of \cite{Tok95}'s $26$-gon. In essence, a polygon $Q$ is taken such that the $45 \mhyphen 45 \mhyphen 90$ triangle $\triangle ABC$ tiles it by successive reflections in a way that the points corresponding to $B$ and $C$ are all vertices of $Q$, and there is a pair of points in the interior of $Q$ that folds down to $A$. Any trajectory on the polygon connecting this pair of points would have to fold down to a trajectory on the triangle that returns to $A$, but such trajectories do not exist (see \cite{Tok95} \textsection 3, Lemma 3.1).

The following lemma extends a lemma by Tokarsky (\cite{Tok95} Lemma 4.1) which is used to construct a more general family of polygons with a pair of points that do not illuminate each other. The above-mentioned result that no trajectory returns to an acute angle in the isosceles triangle can also be derived from this lemma.
\begin{lemma}
\label{lem:tokrep}
Let $Q$ be a triangle $\triangle ABC$ with angles $\angle{A}$ of size $\pi/n$ and $\angle{B}$ of size $m\pi/n$ for some even $n$ and a positive integer $m<n-1$, then there is no billiard trajectory in $Q$ from $A$ coming back to $A$.
\begin{figure}[ht]
\center \includegraphics[width=0.4\textwidth]{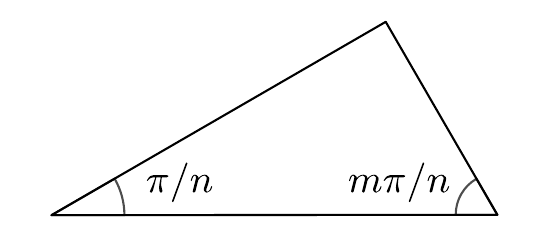}
\caption{}
\label{fig:triangle}
\end{figure}
\end{lemma}
\begin{proof} 
Assume by contradiction that a billiard trajectory $\alpha(t)$ from A to A exists. We may assume that $\alpha(0)=A$ 
and apply a suitable rotation so that $\triangle{ABC}$ is aligned horizontally as in Figure \ref{fig:triangle}. Let $\refl{\theta}$ denote the reflection about a line through the origin which makes an angle of $\theta$ with the $X$-axis and $\rot{\theta}$ the planar rotation by an angle of $\theta$. We can take the dihedral group for the billiard in $\triangle ABC$ to be
$$ \Gamma = \generate{\refl{\frac{\pi}{n}}~,~\refl{0}~,~\refl{\frac{m\pi}{n}}} $$ 
By multiplying the corresponding orthogonal matrices we have $\refl{\frac{\pi}{n}}\refl{0} = \rot{\frac{2\pi}{n}}$ and $(\rot{\frac{2\pi}{n}})^{m-1}\refl{\frac{\pi}{n}}=\refl{\frac{m\pi}{n}}$, so the dihedral group is in fact the symmetry group of the regular $n$-gon:
$$ \Gamma =\generate{\refl{\frac{\pi}{n}}~,~\rot{\frac{2\pi}{n}} } \simeq D_n $$ 
Consider the $2n$ copies of the triangle, unfolded about the vertex $A$ to form a polygon as in Figure \ref{fig:unfolding}. Glued according to the unfolding construction, this is the translation surface associated with $Q$. The vertex $A$ with angle $\pi/n$ unfolds to a removable singularity, considered to be situated at the origin $O$. Let $\widetilde\alpha(t)$ be a lift of $\alpha(t)$ such that for small enough $\eps$, every $0<t<\eps$ has $\alpha(t) \in Q_{(1_\Gamma)}$. Note that at time $t=0$, $\widetilde \alpha(t)$ traverses between $Q_{(1_\Gamma)}$ and the copy of $Q$ rotated by $\pi$, that is $Q_{(\tau)}$ with $\tau = {(\rot{2\pi/n})}^{n/2} \in \Gamma$ being the element of rotation by $\pi$.   
Now clearly for $0 \leq t \leq \eps$ we have $\widetilde\alpha(t) = \tau \widetilde \alpha(-t)$. This symmetry can be extended so that for all $0\leq t \leq T$, $\widetilde{\alpha}(t) = \tau \widetilde{\alpha}(-t)$ with $T$ denoting the minimal period time of $\widetilde\alpha(t)$ as a closed geodesic on a translation surface.
\\
Now observe that for the middle points we have $\widetilde\alpha(T/2) = \widetilde\alpha(-T/2)$ but also $\widetilde\alpha(T/2) = \tau\widetilde\alpha(-T/2)$, so $\widetilde\alpha(T/2)$ is fixed under $\tau$. As $\Gamma$ acts freely on $\inv{\pi}(int(Q))$, $\widetilde\alpha(T/2)$ cannot be a lift of an interior point of $Q$, so it is a lift of a boundary point. By construction, each lift of a point interior to an edge of $Q$ is fixed by a single reflection element of $\Gamma$, so $\widetilde\alpha(T/2)$ must be a singularity (possibly a removable singularity in case the corresponding vertex is $\pi / k$ for some integer $k$). This is a contradiction to $\widetilde\alpha(t)$ being a closed geodesic avoiding singularities (in case $\widetilde\alpha(T/2) \neq O$), or to $T$ being the minimal period (in case $\widetilde\alpha(T/2) = O$).

\begin{figure}[ht]
\center \includegraphics[width=0.7\textwidth]{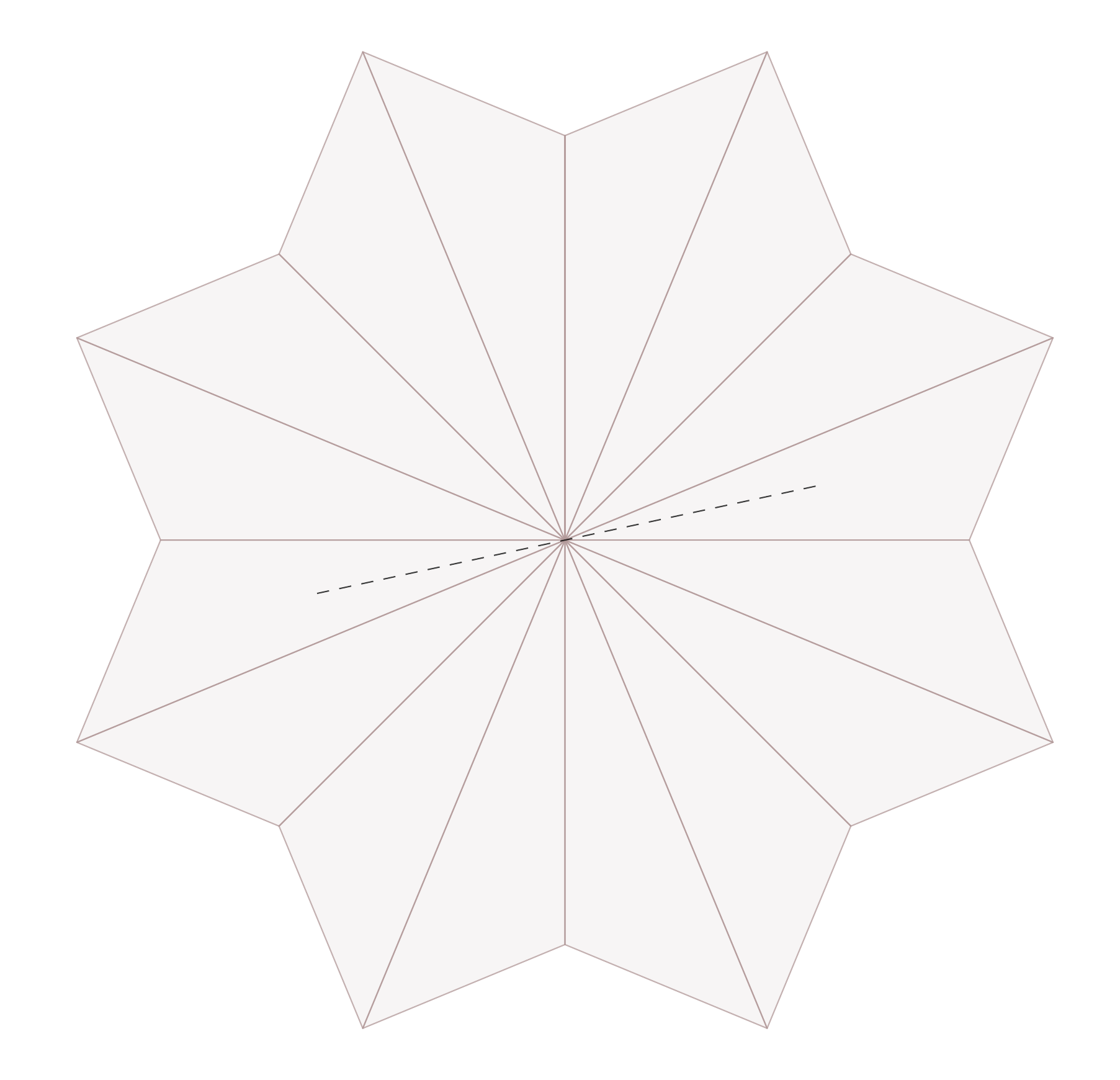}
\caption{An unfolding of the triangle with angles $\pi/8, 5\pi/8$. The dashed line is a short geodesic segment through the removable singularity.}
\label{fig:unfolding}
\end{figure}

\end{proof}

\nocite{*}
\printbibliography
\end{document}